\input amstex
\input epsf
\input prepictex
\input pictex
\input postpictex
\magnification=\magstep1
\documentstyle{amsppt}
\TagsOnRight
\hsize=5.1in                                                  
\vsize=7.8in
\define\R{{\bold R}}
\define\cl{\operatorname{cl}}

\hskip 3 cm {\it Dedicated to S.P. Novikov on the occasion of his 60-th birthday}

\vskip 2cm
\topmatter

\title Lusternik - Schnirelman theory\\
for closed 1-forms\endtitle

\rightheadtext{}
\leftheadtext{}
\author  Michael Farber \endauthor
\address
School of Mathematical Sciences,
Tel-Aviv University,
Ramat-Aviv 69978, Israel
\endaddress
\email farber\@math.tau.ac.il
\endemail
\thanks{The research was supported by a grant from the 
Israel Academy of Sciences and Humanities and by
the Herman Minkowski Center for Geometry} 
\endthanks
\abstract{S. P. Novikov developed an analog of the Morse theory for closed 1-forms.
In this paper we suggest an analog of the Lusternik - Schnirelman theory for closed 1-forms.
For any cohomology class $\xi\in H^1(X,\R)$ we define an integer $\cl(\xi)$
({\it the cup-length associated
with $\xi$}); we prove that any closed 1-form representing $\xi$ has at least $\cl(\xi)-1$ critical
points. The number $\cl(\xi)$ is defined using cup-products in cohomology of some flat line bundles,
such that their monodromy is described by complex numbers, which are not Dirichlet units.}
\endabstract
\endtopmatter

\define\C{{\bold C}}

\define\Z{{\bold Z}}

\define\Hom{\operatorname{Hom}}

\define\Tor{\operatorname{Tor}}

\define\rk{\operatorname{rk}}

\define\coker{\operatorname{coker}}

\define\F{{\Cal F}}


\define\rank{\operatorname{rank}}

\define\E{\Cal E}

\def\<{\langle}
\def\>{\rangle}

\define\pd#1#2{\dfrac{\partial#1}{\partial#2}}

\def\part{\partial}

\NoBlackBoxes

\def\cat{\operatorname {cat}}
\def \cn{\cat_{(N,\partial_+N)}(}
\def \cnm{\cat_{(N_\mu,\partial_+N)}(}
\def\pn{\partial_+N}

\documentstyle{amsppt}

\nopagenumbers

\heading{\bf \S 1. The main result}\endheading

\subheading{1.1} 
Let $X$ be a closed manifold and let $\xi\in H^1(X;\Z)$ be a nonzero cohomology class. 
The Novikov inequalities \cite{N} estimate the numbers of critical points $c_i(\omega)$ of
different indices of any closed 1-form $\omega$ with Morse singularities on $X$ lying in the class $\xi$. 

Novikov type inequalities were constructed in \cite{BF1} for closed 1-forms 
with slightly more general singularities (non-degenerate in the sense of Bott \cite{B}). 
In \cite{BF2}
an equivariant generalization of the Novikov inequalities was found. 

In this paper we will consider the problem of estimating the number of critical points of closed 1-forms
$\omega$ with no non-degeneracy assumption. We suggest here a 
version of the Lusternik - Schnirelman theory for closed 1-forms.

We will define (cf. 1.2 below) a nonnegative integer $\cl(\xi)$, which we will call {\it the cup-length associated with $\xi$}. It is defined in terms of cup-products of some local systems constructed using
$\xi$. 

The main result of the paper consists of the following:

\proclaim{Theorem 1} Let $\omega$ be a closed 1-form on $X$ lying in an integral cohomology class
$\xi\in H^1(X;\Z)$. Let $S(\omega)$ denote the set of critical points of $\omega$, i.e. the set of points $p\in X$
such that $\omega_p=0$. Then the  Lusternik - Schnirelman category of $S(\omega)$ satisfies
$$\cat(S(\omega)) \ge \cl(\xi)-1.\tag1-1$$
In particular, if the set of critical points $S(\omega)$ is finite 
then for the total number $|S(\omega)|$ of the critical points,
$$|S(\omega)|\ge \cl(\xi)-1.\tag1-2$$
\endproclaim

Here $\cat(S)$ denotes the classical Lusternik - Schnirelman category of $S=S(\omega)$, i.e.
the least number $k$, so that $S$ can be covered by $k$ closed subsets 
$A_1\cup\dots \cup A_k$ such that each inclusion $A_j\to S$ is null-homotopic.

A proof of Theorem 1 is given in \S 2.

\subheading{1.2. The cup-length $\cl(\xi)$} Let $\xi\in H^1(X;\Z)$ be an integral cohomology class. 
For any nonzero complex number $a\in \C^\ast$ denote by $\E_a\to X$ the complex flat line bundle
determined by the following condition:
the monodromy along any loop $\gamma\in \pi_1(X)$ is the multiplication by 
$a^{\<\xi, \gamma\>}\in \C$. If $a, b\in \C^\ast$ we have the canonical isomorphism of flat line bundles
$\E_a\otimes \E_b\simeq \E_{ab}$. Therefore we have the cup-product
$$\cup: H^i(X;\E_a)\otimes H^j(X;\E_b)\to H^{i+j}(X;\E_{ab}).\tag1-3$$

\subheading{Definition} {\it The cup-length $\cl(\xi)$ is the largest integer
$k$ such that there exists 
a nontrivial $k$-fold cup product
$$H^{d_1}(X;\E_{a_1})\otimes H^{d_2}(X;\E_{a_2})\otimes \dots \otimes H^{d_k}(X;\E_{a_k})\to 
H^{d}(X;\E_{a}),\tag1-4$$ 
where $d=d_1+\dots+d_k$, $d_1>0, d_2 >0, \dots, d_k>0$, $a=a_1 a_2 \dots a_k$ and among the complex numbers $a_1, \dots , a_k\in \C^\ast$ at least two are not Dirichlet units}.

Recall that {\it a Dirichlet unit} is defined as a complex number $b\ne 0$ such that
$b$ and its inverse $b^{-1}$
are algebraic integers. In other words, Dirichlet units can be characterized as roots of polynomial equations
$$\gamma_0 b^n + \gamma_1 b^{n-1} + \dots + \gamma_{n-1} b + \gamma_n=0,$$
where all $\gamma_i$ are integers and $\gamma_0=\pm 1, \gamma_n=\pm 1.$

Note that the cup-length $\cl(\xi)$ satisfies 
$$0\le \cl(\xi)\le \dim X.\tag1-5$$

In 1.8 we will see examples showing that $\cl(\xi)= \dim X$ is possible.

{}From Theorem 1 we obtain the following simple corollary:

\proclaim{1.3. Corollary} Assume that $X$ is a closed manifold, $\xi\in H^1(X;\Z)$, and 
for some $a\in \C^\ast$, which is not a Dirichlet unit, $H^\ast(X;\E_a)\ne 0$. Then any closed
1-form $\omega$ on $X$ in class $\xi$ has at least one critical point. \endproclaim

\demo{Proof 1} Note that $\cl(\xi)\ge 2$ if $\cl(\xi)>0$. Indeed, assume that
some $H^i(X;\E_a)$ is non-trivial with $a\in \C^\ast$ not a Dirichlet unit. Then by the Poincar\'e
duality there is a non-trivial product $H^i(X;\E_a)\otimes H^{n-i}(X;\E_{a^{-1}})\to H^n(X;\C)$
(where $n=\dim X$),
and so we obtain $\cl(\xi)\ge 2$. Hence by Theorem 1, $|S(\omega)|\ge 1$.\qed\enddemo

\demo{Proof 2}
Corollary 1.3 has also a brief proof independent of Theorem 1, which also explains why 
the definition of the cup-length $\cl(\xi)$ requires that at least two of 
the numbers $a_i$ in (1-4), describing the monodromy, are not Dirichlet units. 
Namely, suppose that 
$\xi\ne 0$ (for $\xi=0$ the form $\omega$ is a function and the statement is trivial) and 
there exists a closed 1-form $\omega$ in class $\xi$ having no critical points. Construct a 
smooth map $f: X\to S^1$ with $\omega=f^\ast(d\theta)$, where $d\theta$ is the angular form on the 
circle. Here 
$$f(x) = \int_{x_0}^x \omega \quad \mod  \Z$$
(we assume that $\xi$ is indivisible). If $\omega$ has no critical points then $f$ is a fibration. Therefore,
$X$ is the mapping torus of a diffeomorphism $h: F\to F$, where $F$ is the fiber of $f$. Hence we obtain (using the Mayer - Vietoris sequence)
$$
\aligned
& H^i(X;\E_a)\simeq \ker [h^\ast -a: H^i(F;\C)\to H^i(F;\C)]\oplus \\
&\coker [h^\ast -a: H^{i-1}(F;\C)\to H^{i-1}(F;\C)]
\endaligned
$$
and now our statement follows from the obvious fact that {\it any eigenvalue of a diffeomorphism 
of a compact manifold, acting on the cohomology, is a Dirichlet unit.} \qed\enddemo

{}From this argument it is clear that Theorem 1 becomes false if we allow the numbers $a_i$ in the
definition of the cup-length in 1.2 to be Dirichlet units. Indeed, one may construct mapping tori
$X$, which admit closed 1-forms with no critical points and may have arbitrarily long non-trivial products (1-4) with $a_i$ Dirichlet units (corresponding to the eigenvalues of the monodromy).
 
\subheading{1.4. Relation to the Novikov numbers}
The following theorem (which is essentially known and stated here only for the sake of completeness)
describes the relation between the cohomology
$H^i(X;\E_a)$ and {\it the Novikov numbers} $b_i(\xi)$, associated with a cohomology class
$\xi\in H^1(X;\Z)$.

\proclaim{Theorem 2} Let $X$ be a closed manifold and let $\xi\in H^1(X;\Z)$ be an integral cohomology class. For fixed $q$ the function
$a\mapsto \dim_\C H^q(X;\E_a),\quad a\in \C^\ast$
has the following behavior:

(a)  it is constant except at finitely many jump points $a = a_1, \dots, a_N$;

(b) the common value of $\dim_\C H^q(X;\E_a)$ for $a\ne a_1, \dots, a_N$ equals the Novikov 
number $b_q(\xi)$;

(c) for $a\in \C^\ast$, being one of the jump points $a_1, \dots, a_N$,
the dimension of the cohomology $\dim_\C H^q(X;\E_a)$ is greater than the Novikov number $b_q(\xi)$;

(d) the jump points $a_1, \dots, a_N$ are algebraic numbers (not necessarily algebraic integers).
\endproclaim

All statements of Theorem 2 except the last one, were announced  in \cite{N3} 
(even in a more general form). 
We will give a simple independent proof in \S 3.

We point out here that as it follows from Theorem 2 {\it for a  transcendental $a\in C$, the dimension of the vector space $H^i(X;\E_a)$ 
does not depend on $a$ and equals the Novikov number $b_i(\xi)$. }

\subheading{1.5. Remarks} 1. A crude estimate for the cup-length $\cl(\xi)$ can be obtained by taking the maximal length of a non-trivial
product (1-4) with $a_1, \dots, a_k$ {\it transcendental.} We will give an example (cf. 1.8, example 3) showing that this
estimate can be really worse than the one provided by Theorem 1. 

2. {\it In the longest nontrivial product (1-4) the number $a$ must be equal 1 and the number $d$ must be equal the dimension of the manifold $\dim X$.}
Indeed, suppose that we have a nontrivial product $v_1\cup\dots\cup v_k\in H^d(X;\E_a)$
with $a=a_1\dots a_k$  not equal 1 or with $d<\dim X,$
and at least two among the numbers $a_1, \dots, a_k$ are not
Dirichlet units. Then $a\ne 1$ implies $d<\dim X$ and (using the Poincar\'e duality) we may find a class
$w\in H^{\dim X-d}(X;\E_{a^{-1}})$ with $v_1\cup\dots\cup v_k\cup w\ne 0$. Hence $\cl(\xi)>k$.
\qed

\subheading{1.6. Forms with non-integral periods} In general, the cohomology class
determined by a closed 1-form $\omega$ belongs to $H^1(X,\R)$, i.e. it has real coefficients. It is clear that multiplying $\omega$ by a non-zero constant $\lambda\ne 0$ does not change the set of critical points $S(\omega)$
and multiplies the cohomology class by $\lambda$. Hence Theorem 1 also gives estimates in the case
of {\it cohomology classes $\xi\in H^1(X,\R)$ of rank 1} (i.e. for classes, which are real multiples of 
integral classes) if we define the associated cup-length $\cl(\xi)$ as follows
$$\cl(\lambda \xi) =\cl(\xi), \quad \lambda\in \R,\,\, \, \lambda \ne 0, \quad \xi \in  H^1(X,\Z).$$ 

Recall, that given a cohomology class $\xi\in H^1(X,\R)$, its {\it rank} is defined as the rank of the abelian group, 
which is the image of the homomorphism $H_1(X,\Z)\to \R$, determined by $\xi$.
Note that the cohomology classes of rank 1 are dense in $H^1(X,\R)$. 
Therefore the following definition makes sense.

\subheading{Definition} Given a class $\xi\in H^1(X,\R)$ of rank $>1$, we define ${\cl}(\xi)$ as the
largest number $k$, such that there exists a sequence of rank 1 classes $\xi_m\in H^1(X,\R)$ with 
$$\cl(\xi_m)\ge k, \qquad \lim_{m\to \infty} \xi_m = \xi,\tag1-7$$ 
and each $\xi_m$, considered as 
a homomorphism $H_1(M;\Z)\to \R$, vanishes on the kernel of the homomorphism
$\xi: H_1(M;\Z)\to \R$.

\proclaim{Theorem 3} Let $\omega$ be a closed 1-form on $X$ lying in a cohomology class
$\xi\in H^1(X;\R)$. Let $S(\omega)$ denote the set of critical points of $\omega$.
Then the  Lusternik - Schnirelman category of $S(\omega)$ satisfies
$$\cat(S(\omega)) \ge {\cl}(\xi)-1.\tag1-8$$
In particular, if the set of critical points $S(\omega)$ is finite 
then for the total number $|S(\omega)|$ of the critical points,
$$|S(\omega)|\ge {\cl}(\xi)-1.\tag1-9$$
\endproclaim

For the proof see \S 3.

\subheading{1.7. Connected sums} Let $X_1$ and $X_2$ be two closed $n$-dimensional manifolds.
We will denote by $X_1\#X_2$ their connected sum. Given cohomology classes
$\xi_\nu\in H^1(X_\nu;\R)$, where $\nu = 1, 2$, the class $\xi_1\# \xi_2\in H^1(X_1\# X_2;\R)$ is well defined (in an obvious way).

In the description of examples (cf. 1.8) we will use the following statement:

\proclaim{Proposition 1} 
$$\cl(\xi_1\# \xi_2) = \max \{ \cl(\xi_1), \cl(\xi_2)\}.\tag1-10$$
\endproclaim

Proof is given in \S 3.

\subheading{1.8. Examples} 1. In the notations of the previous subsection, let $\xi_1=0$ and 
suppose that $\xi_2\ne 0$ can be realized by a closed 1-from with no critical points (for example,
fibration over the circle). Then we obtain from Proposition 1
that $\cl(\xi_1\# \xi_2)=\cl(\xi_1)$. Since $\xi_1=0$, the cup-length $\cl(\xi_1)$ is the
usual cup-length of the manifold $X_1$ with rational coefficients. 

To have a specific example, let us take $X_1=T^n$, $X_2=S^1\times S^{n-1}$, $\xi_1=0$ and 
$\xi_2\in H^1(X_2;\Z)$ being a generator. Then we have for $\xi=\xi_1\# \xi_2\in H^1(X_1\# X_2;\R)$
$$\cl(\xi_1\# \xi_2)=n.\tag1-11$$
Therefore, by Theorem 1, any closed 1-form $\omega$ on $X_1\# X_2$ lying in class 
$\xi$ has a least
$n-1$ critical points.

2. In a similar way one may construct examples of cohomology classes of higher rank with many 
critical points. Namely, suppose that $X_1=T^n$ and $\xi_1=0$; take for $X_2$ arbitrary closed
manifold of dimension $n$ with a cohomology class $\xi_2\in H^1(X_2;\R)$ of rank $q$.
Then for the class $\xi=\xi_1\# \xi_2\in H^1(X_1\# X_2;\R)$ 
(having rank $q$) we again obtain  
$\cl(\xi)=n$ (by comparing Proposition 1 with (1-5).

One may take, for example, $X_2=T^q\times S^{n-q}$ with $\xi_2$ induced
from a maximally irrational class on the torus $T^q$.

3. Let $X$ be a 3-dimensional manifold obtained by 0-framed surgery on the knot $5_2$:


%
$$
\beginpicture
\setcoordinatesystem units <1.00000cm,1.00000cm>
\linethickness=1pt
\setlinear
\linethickness= 0.500pt
%
%
%
\plot	 8.541 22.320  8.779 22.368
 	 8.895 22.383
	 9.005 22.380
	 9.109 22.359
	 9.207 22.320
	 9.301 22.272
	 9.390 22.221
	 9.475 22.168
	 9.557 22.114
	 9.627 22.052
	 9.680 21.979
	 9.715 21.894
	 9.731 21.796
	 9.731 21.686
	 9.715 21.562
	 9.702 21.495
	 9.684 21.424
	 9.662 21.350
	 9.636 21.273
	 9.607 21.193
	 9.575 21.114
	 9.539 21.034
	 9.501 20.955
	 9.460 20.876
	 9.416 20.796
	 9.369 20.717
	 9.319 20.637
	 9.269 20.562
	 9.222 20.493
	 9.140 20.376
	 9.072 20.286
	 9.017 20.225
	 8.961 20.169
	 8.890 20.098
	 8.803 20.010
	 8.700 19.907
	 /
\plot  8.700 19.907  8.477 19.685 /
\linethickness= 0.500pt
%
%
%
\plot	 7.747 19.685  8.001 19.526
 	 8.068 19.486
	 8.142 19.444
	 8.223 19.400
	 8.311 19.356
	 8.405 19.309
	 8.507 19.261
	 8.616 19.212
	 8.731 19.161
	 8.846 19.114
	 8.951 19.074
	 9.046 19.042
	 9.132 19.018
	 9.209 19.002
	 9.276 18.994
	 9.382 19.002
	 9.459 19.042
	 9.513 19.114
	 9.546 19.217
	 9.554 19.280
	 9.557 19.352
	 9.554 19.427
	 9.545 19.500
	 9.530 19.573
	 9.509 19.645
	 9.482 19.716
	 9.450 19.786
	 9.411 19.855
	 9.366 19.923
	 /
\plot  9.366 19.923  9.176 20.193 /
\linethickness= 0.500pt
%
%
%
\plot	 7.906 20.542  7.763 20.431
 	 7.680 20.363
	 7.572 20.268
	 7.510 20.212
	 7.441 20.148
	 7.367 20.079
	 7.287 20.003
	 7.207 19.927
	 7.135 19.860
	 7.013 19.748
	 6.921 19.669
	 6.858 19.622
	 6.779 19.514
	 6.751 19.428
	 6.731 19.320
	 6.727 19.208
	 6.747 19.110
	 6.791 19.025
	 6.858 18.955
	 6.942 18.902
	 7.037 18.871
	 7.141 18.862
	 7.255 18.875
	 7.373 18.901
	 7.489 18.931
	 7.603 18.965
	 7.715 19.002
	 /
\plot  7.715 19.002  7.938 19.082 /
\linethickness= 0.500pt
%
%
%
\plot	 6.985 20.320  6.763 20.542
 	 6.664 20.662
	 6.623 20.729
	 6.588 20.800
	 6.559 20.876
	 6.537 20.956
	 6.520 21.041
	 6.509 21.130
	 6.503 21.221
	 6.501 21.312
	 6.503 21.403
	 6.509 21.495
	 6.519 21.586
	 6.533 21.677
	 6.550 21.769
	 6.572 21.860
	 6.602 21.948
	 6.645 22.030
	 6.699 22.105
	 6.767 22.173
	 6.846 22.236
	 6.938 22.291
	 7.043 22.341
	 7.160 22.384
	 7.280 22.420
	 7.395 22.449
	 7.504 22.472
	 7.608 22.487
	 7.707 22.495
	 7.800 22.497
	 7.887 22.491
	 7.969 22.479
	 8.045 22.463
	 8.113 22.445
	 8.227 22.408
	 8.312 22.366
	 8.366 22.320
	 8.411 22.260
	 8.465 22.173
	 8.530 22.061
	 8.566 21.996
	 8.604 21.923
	 8.638 21.848
	 8.662 21.775
	 8.674 21.702
	 8.676 21.630
	 8.666 21.559
	 8.646 21.489
	 8.615 21.420
	 8.572 21.352
	 /
\plot  8.572 21.352  8.382 21.082 /
\linethickness= 0.500pt
%
%
%
\plot	 8.128 22.098  7.938 21.860
 	 7.863 21.734
	 7.842 21.666
	 7.830 21.594
	 7.830 21.519
	 7.839 21.440
	 7.859 21.358
	 7.890 21.273
	 7.925 21.189
	 7.959 21.114
	 7.993 21.046
	 8.025 20.987
	 8.086 20.892
	 8.144 20.828
	 8.206 20.779
	 8.283 20.729
	 8.373 20.676
	 8.477 20.622
	 /
\plot  8.477 20.622  8.700 20.510 /
\linethickness=0pt
\putrectangle corners at  6.445 22.543 and  9.779 18.796
\endpicture
$$
\centerline{Figure 1.}
\smallskip

This knot has Alexander polynomial $\Delta(\tau) = 2-3\tau +2\tau^2$.
Then $H^1(X;\Z)=\Z$ and taking $\xi\in H^1(X;\Z)$ to be a generator we find that $H^1(X;\E_a)$
is trivial for all $a\in \C^\ast$, which are not the roots of the Alexander polynomial. It is easy to check
that if $a$ is one of the 
roots of $2-3a+2a^2=0$ then $H^1(X;\E_a)\ne 0$. Note that the roots of $2-3a+2a^2=0$ are not Dirichlet units.
Hence we obtain (using Theorem 2) that all Novikov
numbers are trivial, however by Corollary 1.3 we obtain that any closed 1-forms in the 
class $\xi$ has at least 1 critical point. 

4. Let $X_g$ be a compact Riemann surface of genus $g>1$. Then for any $\xi\in H^1(X_g;\Z)$
and $a\ne 1$ holds $\dim H^1(X_g;\E_a)=2g-2$ and all other cohomology groups are trivial. Then we
obtain that $\cl(\xi)=2$ (cf. Proof 1 of Corollary 1.3).  Thus Theorem 1 predicts existence of 
one critical point of any closed 1-form in any non-trivial cohomology class (which also follows from 
the Hopf's theorem). 

We observe that {\it 
for any $\xi\in H^1(X_g;\Z)$ with $\xi\ne 0$, $g>1$ there exists a closed 1-form on the surface $X_g$ lying in class
$\xi$ and having precisely one critical point.} 

Indeed, it is known that on the surface
$X_{g-1}$ of genus $g-1$ there exists a smooth 
function with precisely three critical points: a maximum, a minimum and a saddle point, cf. \cite{DNF},
chapter 2, figure 81.
Let $f:X_{g-1}\to \R$ be such function. We may assume that the critical values of $f$ are
$-1, 0, 1$. The level sets $f^{-1}(-1/2)=C_{-1}$ and $f^{-1}(1/2)=C_1$ are circles. If one picks an orientation of the surface $X_{g-1}$, these circles become oriented; if the orientation of $X_{g-1}$
is reversed then the orientations of $C_{-1}$ and $C_1$ are reversed. 
Let now $X_g$
be obtained from $f^{-1}([-1/2,1/2])$ by identifying the points of the circles $C_{-1}$ and $C_1$
by a diffeomorphism $C_{-1}\to C_1$ preserving the orientations (the meaning of this is clear from the 
remark above).
We obtain a map $h:X_g\to S^1=\R/\Z$ by $h(x)=f(x)\,  \mod\, \Z$ for $x\in X_g$. This gives a
map to the circle $S^1$ with precisely one critical point. Clearly, the closed 1-form 
$h^\ast(d\theta)$ has only one critical point.

Hence our estimates are exact in the case of surfaces.

\heading{\bf \S 2. Proof of Theorem 1}\endheading

\subheading{2.1}
Since we assume that the cohomology class $\xi$ of $\omega$ is integral, $\xi \in H^1(X,\Z)$,
there is a smooth map $f: X \to S^1$ such that $\omega =f^\ast(d\theta)$, where $d\theta$ is the 
standard angular form on the circle $S^1\subset \C$, $S^1=\{z; |z|=1\}.$ 
We will also assume (without loss of generality) that $\xi$ is indivisible and that $1\in S^1$ is a regular
value of $f$.

Denote $f^{-1}(1)$ by $V\subset X$; it is a smooth codimension 
one submanifold.  
Let $N$ denote the manifold obtained by cutting of $X$ along $V$. We will denote by $\partial_+N$
and $\partial_-N$ the components of the boundary of $N$. We get a smooth function 
$$g: N \to [0,1],\tag2-1$$
so that 
\roster
\item $g$ is identically 0 on $\partial_+N$ and identically 1 on $\partial_-N$ and has no critical points
near $\partial N$;
\item letting $p: N\to X$ denote the natural identification map, then for any $x\in N$,
$$f(p(x)) = \exp(2\pi i g(x));$$
\item $p$ maps the set $S(g)$ of critical points of $g$ homeomorphically onto the set $S(\omega)$
of critical points of $\omega$.
\endroster

$$
\beginpicture
\setcoordinatesystem units <1.00000cm,1.00000cm>
\linethickness=1pt
\setlinear
%
%
\linethickness= 0.500pt
\putrule from  2.508 22.892 to  2.508 19.717
%
%
\linethickness= 0.500pt
\putrule from  7.588 20.987 to  8.858 20.987
%
%
\plot  8.604 20.923  8.858 20.987  8.604 21.050 /
%
%
%
\linethickness= 0.500pt
\putrule from 10.763 19.558 to 10.763 24.162
%
%
\plot 10.827 23.908 10.763 24.162 10.700 23.908 /
%
%
%
\linethickness= 0.500pt
\setdashes < 0.1270cm>
\plot  4.413 19.717 10.763 19.717 /
%
%
\linethickness= 0.500pt
\setsolid
\putrule from 10.604 19.717 to 10.922 19.717
%
%
\linethickness= 0.500pt
\plot 10.604 23.527 10.604 23.527 /
%
%
\linethickness= 0.500pt
\putrule from 10.604 22.892 to 10.922 22.892
%
%
\linethickness= 0.500pt
\setdashes < 0.1270cm>
\plot  4.413 22.892 10.763 22.892 /
\linethickness= 0.500pt
\setsolid
%
%
\plot  2.508 22.892 	 2.573 22.880
	 2.650 22.862
	 2.736 22.838
	 2.826 22.812
	 2.916 22.786
	 3.002 22.763
	 3.079 22.744
	 3.143 22.733
	 3.207 22.727
	 3.284 22.723
	 3.370 22.720
	 3.461 22.719
	 3.552 22.720
	 3.638 22.723
	 3.715 22.727
	 3.778 22.733
	 3.849 22.742
	 3.942 22.756
	 4.046 22.773
	 4.154 22.793
	 4.254 22.816
	 4.337 22.840
	 4.413 22.892
	 4.375 22.933
	 4.274 22.968
	 4.206 22.983
	 4.129 22.996
	 4.044 23.008
	 3.955 23.018
	 3.862 23.027
	 3.768 23.035
	 3.676 23.041
	 3.586 23.045
	 3.502 23.049
	 3.425 23.051
	 3.302 23.050
	 3.213 23.044
	 3.097 23.030
	 3.033 23.020
	 2.966 23.009
	 2.899 22.998
	 2.832 22.986
	 2.707 22.960
	 2.604 22.934
	 2.508 22.892
	 2.508 22.892
	/
\linethickness= 0.500pt
%
%
\plot  2.508 19.717 	 2.635 19.680
	 2.745 19.650
	 2.839 19.624
	 2.920 19.603
	 2.988 19.587
	 3.047 19.574
	 3.143 19.558
	 3.207 19.552
	 3.284 19.548
	 3.370 19.545
	 3.461 19.544
	 3.552 19.545
	 3.638 19.548
	 3.715 19.552
	 3.778 19.558
	 3.874 19.574
	 3.933 19.587
	 4.002 19.603
	 4.082 19.624
	 4.177 19.650
	 4.286 19.680
	 4.413 19.717
	/
\linethickness= 0.500pt
\setdashes < 0.1270cm>
%
%
\plot  2.508 19.717 	 2.607 19.738
	 2.698 19.758
	 2.784 19.776
	 2.863 19.792
	 2.936 19.806
	 3.004 19.819
	 3.125 19.841
	 3.228 19.857
	 3.317 19.867
	 3.393 19.874
	 3.461 19.876
	 3.528 19.874
	 3.605 19.867
	 3.694 19.857
	 3.797 19.841
	 3.918 19.819
	 3.986 19.806
	 4.059 19.792
	 4.138 19.776
	 4.223 19.758
	 4.315 19.738
	 4.413 19.717
	/
\linethickness= 0.500pt
\setsolid
%
%
\plot  4.445 22.860 	 4.443 22.781
	 4.442 22.707
	 4.441 22.638
	 4.441 22.575
	 4.441 22.461
	 4.444 22.364
	 4.449 22.282
	 4.456 22.211
	 4.477 22.098
	 4.510 21.992
	 4.562 21.866
	 4.627 21.744
	 4.699 21.654
	 4.807 21.589
	 4.876 21.562
	 4.949 21.538
	 5.023 21.516
	 5.094 21.497
	 5.207 21.463
	 5.313 21.435
	 5.378 21.419
	 5.447 21.403
	 5.516 21.383
	 5.581 21.361
	 5.683 21.304
	 5.736 21.240
	 5.782 21.153
	 5.819 21.063
	 5.842 20.987
	 5.854 20.917
	 5.861 20.830
	 5.858 20.742
	 5.842 20.669
	 5.789 20.579
	 5.708 20.486
	 5.615 20.405
	 5.524 20.352
	 5.493 20.352
	 5.415 20.342
	 5.317 20.327
	 5.220 20.308
	 5.143 20.288
	 5.034 20.249
	 4.968 20.223
	 4.899 20.194
	 4.831 20.163
	 4.767 20.130
	 4.667 20.066
	 4.566 19.952
	 4.499 19.853
	 4.459 19.790
	 4.413 19.717
	/
\linethickness= 0.500pt
%
%
\plot  4.159 20.860 	 4.234 20.809
	 4.299 20.767
	 4.408 20.708
	 4.496 20.677
	 4.572 20.669
	 4.657 20.685
	 4.752 20.728
	 4.869 20.806
	 4.938 20.859
	 5.017 20.923
	/
\linethickness= 0.500pt
%
%
\plot  4.255 20.828 	 4.348 20.858
	 4.417 20.877
	 4.508 20.892
	 4.611 20.878
	 4.689 20.858
	 4.794 20.828
	/
%
%
\put{$0$} [lB] at 11.144 19.590
%
%
%
%
%
%
%
%
\put{$g$} [lB] at  8.223 21.431
%
%
%
%
%
%
%
%
%
%
\put{$\partial_+N$}[lB] at 1.143 19.463
\put{$\partial_-N$}[lB] at 1.270 22.892
\put{Figure 2} [lB] at  5.683 18.447
%
%
\put{$1$} [lB] at 11.144 22.765
\linethickness=0pt
\putrectangle corners at  1.143 24.225 and 11.144 18.383
\endpicture
$$
\bigskip

\subheading{2.2}
For any subset $X\subset N$ containing $\pn$ 
we will denote by $\cn X)$ the minimal number $k$ such that $X$
can be covered by $k+1$ closed subsets 
$$X\subset A_0\cup A_1\cup A_2\cup \dots \cup A_k\subset N$$
with the following properties:
\roster
\item $A_0$ is {\it a collar of $\pn$}, i.e. it contains a neighborhood of $\pn$ and 
the inclusion $A_0\to N$ is homotopic to a map $A_0\to \pn$ keeping the points of $\pn\subset A$
fixed;
\item for $j=1, 2, \dots, k,$ each set $A_j$ is disjoint from $\pn$ and the inclusion $A_j\to N$
is null-homotopic.
\endroster

 The number $\cn X)$ can be viewed as a relative version of the Lusternik
- Schnirelman category.

Our purpose in this subsection is to prove the inequality
$$\cat S(\omega) \ge \cn N).\tag2-2$$
The arguments here are modifications of the standard arguments. 

We will need Lemmas 1 - 4.

\proclaim{Lemma 1} Let $X, Y\subset N$ be two closed subsets containing a neighborhood of
$\pn$. Suppose that there exists a deformation $G_t: Y\to N$, $t\in [0,1]$, such that $G_0=\text{inclusion}:Y\to N$, $G_1(Y)\subset X$
and $G_t(x)=x$ for all $x\in \partial_+N$, $t\in [0,1]$. Then 
$$\cn Y) \le \cn X).\tag2-3$$
\endproclaim
\demo{Proof} Suppose that $\cn X)=k$ and let $A_0\cup A_1\cup \dots \cup A_k$ be a cover of $X$ by closed subsets as above.
Set $B_0=G_1^{-1}(A_0)$ and for $j=1,2, \dots, k$ let $B_j$ be defined as 
$G_1^{-1}(A_j)$, with a small cylindrical neighborhood of $\pn$ removed. Let us show that the sets $B_j$, $j=0, 1,2, \dots, k$ satisfy the requirements of the above definition. 
We have the following deformation
$G_t|_{B_j}: B_j\to N$, $t\in [0,1]$, which starts with the inclusion $B_j\to N$ and ends with
a map $B_j\to A_j$. After that for $j>0$ we may apply the deformation which shrinks $A_j$ to a point; for $j=0$ we apply the deformation which brings $A_0$ to $\pn$ keeping 
$ \pn\subset A_0$ fixed.
This gives a covering of $Y$ with the required properties. Therefore, $\cn Y)\le k$.
 \qed
\enddemo
\proclaim{Lemma 2} (a) Let $A\subset N-\pn$ be a compact subset such that the inclusion $A\to N$
is null-homotopic. Then for any $\epsilon>0$ small enough
the $\epsilon$-neighborhood $U_\epsilon(A)$ of $A$ is also null-homotopic in $N$. 
(b) Let $A$ be a closed subset containing a neighborhood of $\pn$, such that $A$ can be deformed into $\pn$ in $N$ keeping the points of $\pn$ fixed. Then for small $\epsilon$ the 
$\epsilon$-neighborhood 
$U_\epsilon(A)$ of $A$ can also be deformed into $\pn$ keeping the points of $\pn$ fixed.
\endproclaim
\demo{Proof} Let us prove (a); the proof for (b) is the same. 
Suppose that $G: A\times [0,1]\to N$ is the given deformation of $A$ into $\pn$. 
Embed $N$ into some Euclidean space $\R^n$. $G$ defines a mapping from the closed subset
$A\times [0,1]\subset N\times [0,1]$ to $\R^n$. By the Tietze theorem there exists a continuous
map $\tilde G: N\times [0,1]\to \R^n$ which coincides with $G$ on $A\times [0,1]$. Let
$U$ be a small tubular neighborhood of $N$ in $\R^n$ and let $\pi: U\to N$ be the projection (retraction). Suppose that $\epsilon$ is so small that the $\epsilon$-neighborhood of 
$A\times [0,1]$ is contained
in $\tilde G^{-1}(U)$. Then the image $\tilde G(U_\epsilon(A\times [0,1]))$ is contained in $U$ and, composing $\tilde G$ with the projection $\pi$, we obtain a deformation
$U_\epsilon(A)\times [0,1] \to N$. The final map of this deformation brings $A$ into $\pn$ and so 
the rest of $U_\epsilon(A)$ is contained in a collar of $\pn$. Hence  
$\tilde G(U_\epsilon(A\times 1))$ could be brought into $\pn$
by another deformation.\qed

\enddemo

$$\epsffile{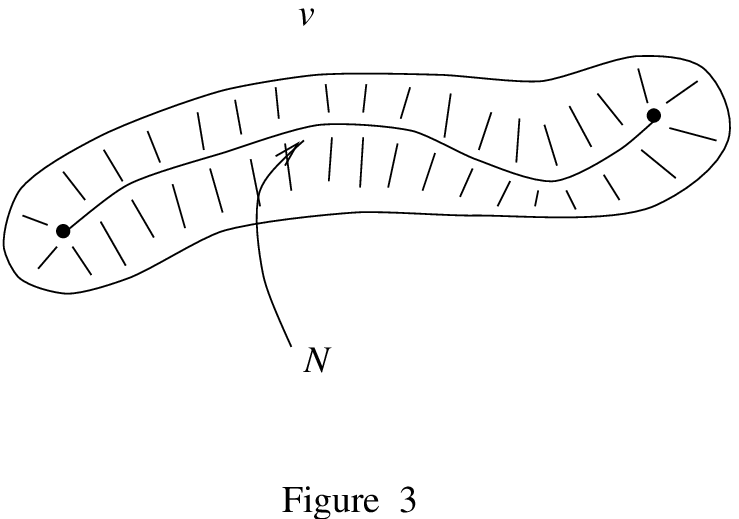}$$

\proclaim{Lemma 3} For any closed subset $X\subset N$ containing $\pn$ 
there exists $\epsilon >0$ such that
$$\cn U_\epsilon (X))= \cn X).\tag2-4$$
\endproclaim
\demo{Proof} The inequality $\cn U_\epsilon (X))\ge \cn X)$
is obvious. The opposite inequality follows from Lemma 2.\qed
\enddemo
\proclaim{Lemma 4} For a pair of closed subsets $X, Y\subset N$ such that $X$ contains
$\pn$, it holds
$$\cn X) + \cat_N(Y) \ge \cn X\cup Y).\tag2-5$$
\endproclaim
Here $\cat_N(Y)$ denotes the Lusternik-Schnirelman category of $Y$ with respect to $N$, i.e. the 
minimal numbers $k$ such that $Y$ can be covered by $k$ closed subsets which are null-homotopic
in $N$.
\demo{Proof} Obvious. \qed\enddemo

We are in a position now {\it to prove the inequality} (2-2). 

We may assume that $S(\omega)= S(g)$ has only finitely many connected components, 
since otherwise 
$\cat(S(\omega))$ is infinite, and Theorem 1 is obviously true. Also, the function $g$ (defined in 2.1)
is constant on
each connected component of the critical point set $S(g)$ 
(this follows from the Sard's theorem since the 
image of a connected component under $g$ must be connected and must have measure zero).
Hence we obtain that the set of critical values of function $g$ is finite.
 
Consider the function
$$F(\mu) = \cnm N_\mu),\quad\text{where}\quad N_\mu= g^{-1}([0,\mu]),\qquad \mu \in [0,1].$$
\midinsert
$$\epsffile{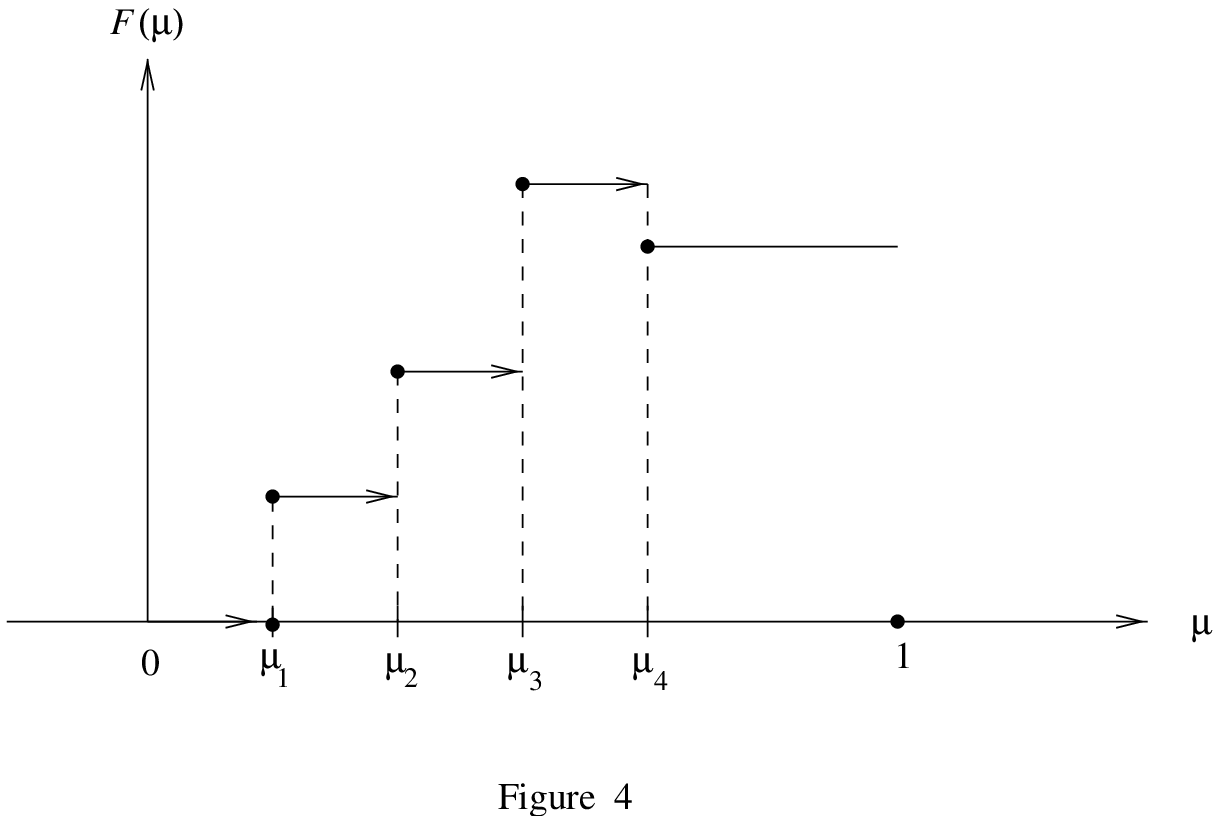}$$
\endinsert

\noindent
It is clear that $F(\mu) =0$ for $\mu\ge 0$ small, and $F(\mu) =\cn N)$ for $\mu\le 1$ close to 1.
Moreover, from the basic theorems of the Morse theory it follows that $F(\mu)$ is constant on each 
segment $[\mu,\mu']$, containing no critical values. Thus, $F(\mu)$ is a step function which
may have finitely many jumps $\mu_1\le \mu_2\le \dots\le \mu_s$ and all the jump points
$\mu_1\le \mu_2\le \dots\le \mu_s$ are the critical values of the function $g$. 

We want to show that for any critical value $\mu\in [0,1]$,
$$F(\mu) - F(\mu -\epsilon)\le \cat (S_\mu),\tag2-6$$
where $S_\mu$ is the set of all critical points of $g$ on the level $g^{-1}(\mu)$. 
Here $\epsilon>0$ is so small that $[\mu-\epsilon, \mu) $ contains no critical values.
To prove (2-6) we observe that choosing $\delta>0$ small enough we get (by Lemma 4)
$$\cnm N_\mu-U_\delta(S_\mu)) + \cat_{N_\mu}(\overline U_\delta(S_\mu)) \ge F(\mu)\tag2-7$$
and also by Lemma 1 (using the deformation determined by the gradient flow of the function $g$)
$$
\aligned
&\cnm N_\mu-U_\delta(S_\mu)) \le \cnm N_{\mu-\epsilon}) \le\\
& \cat_{(N_{\mu-\epsilon}, \partial_+N)} N_{\mu-\epsilon} = F(\mu -\epsilon).
\endaligned\tag2-8$$
In addition we have 
$$\cat_{N_\mu}(\overline U_\delta(S_\mu)) = \cat_{N_\mu}(S_\mu) \le \cat(S_\mu).\tag2-9$$
Combining (2-7), (2-8) and (2-9) proves (2-6).

The total jump of the function $F$ on the interval $[0,1]$ 
equals $\cn N)$ and by (2-6) at each critical value $\mu$ 
the value of the jump does not
exceed the Lusternik - Schnirelman category $\cat(S_\mu)$ of the set of critical points on the level.
(Note that $F(\mu)$ may also have some negative jumps.) 
Hence the total jump $\cn N)$ does not exceed 
$$\sum_\mu  \cat(S_\mu) = \cat(S(g))=\cat (S(\omega)).\tag2-10$$
This proves (2-2).

\subheading{2.3. The deformation complex} Later (cf. 2.7) we will prove that
$$\cn N) \ge \cl(\xi) - 1.\tag2-11$$
Together with (2-2) this will complete the proof of the Theorem.
The proof of (2-11) will consist of building a {\it polynomial deformation} of the cochain complex
$C^\ast(X;\E_a)$ (where $a$ is viewed as a parameter) into $C^\ast(N,\pn)$ as $a\to \infty$. 
The {\it deformation} understood here as a finitely generated free cochain complex  $C^\ast$
over the ring $P=\Z[\tau]$ of polynomials with integral coefficients satisfying (a) and (b) below.

The construction of the deformation $C^\ast$ goes as follows. 
We shall assume that $N$ is triangulated
and $\partial N$ is a subcomplex. Recall that $V=f^{-1}(1)$ (cf. {2.1}) and we will denote 
by $i_\pm : V \to N$ the inclusions, which identify $V$ with $\partial_\pm N$ correspondingly.
Denote by $C^q(N)$ and $C^q(V)$ the free abelian groups of integer valued cochains
and $\delta_N: C^q(N) \to C^{q+1}(N)$ and by $\delta_V: C^q(V) \to C^{q+1}(V)$  
the coboundary homomorphisms. 

Let $C^q(N)[\tau]$ and $C^{q-1}(V)[\tau]$ denote free $P$-modules formed by "polynomials
with coefficients" in the corresponding abelian groups; for example, an element $c\in  C^q(N)[\tau]$
is a formal sum $c=\sum_{i\ge 0}c_i\tau^i$ with $c_i\in C^q(N)$ and only finitely many $c_i$'s
are nonzero. We shall consider of $c\in C^q(N)[\tau]$ as a polynomial (complex) curve,
which associates a point in $C^q(N)$ with a complex number $\tau \in \C$. The $P$-module structure
is given naturally as follows: $\tau\cdot c= \sum_{i\ge 0}c_i\tau^{i+1}$. It is clear that 
$C^q(N)[\tau]$ and $C^{q-1}(V)[\tau]$ are free finitely generated $P$-modules and their ranks
equal to the number of $q$-dimensional simplices in $N$ or $(q-1)$-dimensional simplices
in $V$, correspondingly.

Consider the natural $P$-module
extensions
$$\delta_N: C^q(N)[\tau] \to C^{q+1}(N)[\tau],\quad\text{ and}\quad 
\delta_V: C^q(V)[\tau] \to C^{q+1}(V)[\tau].\tag2-12$$
They act coefficientwise so that $\delta_N$ and $\delta_V$ are $P$-homomorphisms. 
For example, if $\alpha=\sum_{i\ge 0}\alpha_i\tau^i\in C^q(N)[\tau]$ then 
$\delta_N(\alpha) = \sum_{i\ge 0}\delta_N(\alpha_i)\tau^i$.

Now we define a finitely generated free cochain complex  $C^\ast$
over the ring $P=\Z[\tau]$ of polynomials with integral coefficients as follows:
$C^\ast=\oplus C^q$, where
$$C^q = C^q(N)[\tau]\oplus C^{q-1}(V)[\tau].\tag2-13$$
Elements of $C^q$ will be denoted as pairs $(\alpha,\beta)$, where 
$\alpha\in C^q(N)[\tau]$ and $\beta\in C^{q-1}(V)[\tau]$.
The differential
$\delta: C^q\to C^{q+1}$ is given by the following formula:
$$\delta(\alpha, \beta) = (\delta_N(\alpha), (i_+^\ast -\tau i_-^\ast)(\alpha) -\delta_V(\beta)),\tag2-14$$
where $\alpha\in C^q(N)[\tau]$ and $\beta\in C^{q-1}(V)[\tau]$. It is clear that $C^\ast$ is the usual
cylinder of the chain map $i^\ast_+-\tau i^\ast_-$ with a shifted grading.

We claim now that:{\it 

(a) for any nonzero complex number $a\in \C^\ast$ there is a canonical isomorphism}
$$
\CD
E^q_a: H^q(C^\ast\otimes _P\C_a) \overset{\simeq}\to\longrightarrow H^q(X;\E_{a^{-1}}).
\endCD
\tag2-15$$
Here $\C_a$ is $\C$ which is viewed as a $P$-module with the following structure:
$\tau x= ax$ for $x\in \C$.
We will call $E^q_a$ {\it an isomomorphism of evaluation at $\tau=a$};

(b) {\it for $a=0$ we also have a canonical evaluation isomorphism}
$$E^q_0: H^q(C^\ast\otimes _P\Z_0) \to H^q(N, \partial_+N;\Z),\tag2-16$$
where $\Z_0$ is $\Z$ with the following $P$-module structure: $\tau x=0$ for any $x\in \Z$.

To show (a) we note that $H^q(X;\E_{a^{-1}})$ can be identified with the cohomology of 
complex $C^\ast(X;\E_{a^{-1}})$ consisting of cochains $\alpha\in C^q(N)$ satisfying the boundary
conditions
$$i_-^\ast(\alpha) = a^{-1}i_+^\ast(\alpha) \in C^q(V).$$
The complex $C^\ast\otimes_P\C_a$ can be viewed as
$$C^q\otimes_P\C_a = C^q(N)\oplus C^{q-1}(V)$$ with the differential given by
$$\delta(\alpha, \beta) = (\delta_N(\alpha), (i_+^\ast -a i_-^\ast)(\alpha) -\delta_V(\beta)),\tag2-17$$
where $\alpha\in C^q(N)$ and $\beta\in C^{q-1}(V)$. It is clear that there is a chain homomorphism
$C^\ast(X;\E_{a^{-1}})\to C^\ast\otimes_P\C_a$ (acting by $\alpha\mapsto (\alpha,0))$. 
It is easy to see that it induces an isomorphism on the cohomology. Indeed, suppose that a cocycle
$\alpha\in C^q(X;\E_{a^{-1}})$ bounds in the complex $C^\ast\otimes_P\C_a$ then there are
$\alpha_1\in C^{q-1}(N)$, $\beta_1\in C^{q-2}(V)$ such that 
$\alpha=\delta_N(\alpha_1)$, $i_+^\ast(\alpha_1) - ai_-^\ast(\alpha_1)-\delta_V(\beta_1)=0.$
We may find a cochain $\beta_2\in C^{q-2}(N)$ such that $i_+^\ast(\beta_2)=\beta_1$ and 
$i_-^\ast(\beta_2)=0$ (by extending $\beta_1$ into a neighborhood of $\pn$). Then setting
$\alpha_2=\alpha_1-\delta_N(\beta_2)$ we have 
$$\alpha=\delta_N(\alpha_2),\qquad i_+^\ast(\alpha_2) - ai_-^\ast(\alpha_2)=0,\tag2-18$$
which means that $\alpha$ also bounds in $C^q(X;\E_{a^{-1}})$. 

Similarly, suppose that $(\alpha, \beta)$ is a cocycle of complex $C^\ast\otimes_P\C_a$. As above we
may find a cochain $\beta'\in C^{q-1}(N)$ with $i_+^\ast(\beta')=\beta$ and 
$i_-^\ast(\beta')=0$. Then $(\alpha-\delta_N(\beta'), 0)$ is cohomologous to the initial cocycle 
$(\alpha, \beta)$ and it is a cocycle of $C^\ast(X;\E_{a^{-1}})$.

This proves (a). The statement  (b) follows similarly.\qed

\subheading{2.4. Relative deformation complex} We will define now a relative version of the deformation complex $C^\ast$.

Let $A\subset N$ be a simplicial subcomplex. We will assume that $A$ {\it is disjoint from} $\partial_+N$. Let $C^q(N,A)$ denote the free abelian group of integer-valued cochains on 
$N$ which vanish on $A$. Let $C^q(N, A)[\tau]$ be constructed similarly to 
$C^q(N)[\tau]$, cf. above. We define the complex $C^\ast_A$ as follows:
$$C^q_A = C^q(N,A)[\tau]\oplus C^{q-1}(V)[\tau].\tag2-19$$
The differential
$\delta: C^q_A\to C^{q+1}_A$ is defined by the following formula:
$$\delta(\alpha, \beta) = 
(\delta_{N,A}(\alpha), (i_+^\ast -\tau i_-^\ast)(\alpha) -\delta_V(\beta)),\tag2-20$$
where $\alpha\in C^q(N,A)[\tau]$ and $\beta\in C^{q-1}(V)[\tau]$. 
Here
$\delta_{N,A}: C^q(N,A) \to C^{q+1}(N,A)$ and $\delta_V: C^q(V) \to C^{q+1}(V)$  denote
the coboundary homomorphisms and also their $P$-module extension. 
$i_\pm^\ast : C^q(N,A) \to C^q(V)$
denote the restriction maps of chains, and the same symbols denote also their polynomial extensions
$i_\pm^\ast : C^q(N,A)[\tau] \to C^q(V)[\tau]$. 

Similarly to statements (a) and (b) in 2.3 we have:

{\it (a') for any $a\in \C^\ast$ there is a natural isomorphism
$$H^i(C^\ast_A\otimes_P \C_a) \simeq H^i(X,p(A);\E_{a^{-1}}),\tag2-21$$
where $p:N \to X$ is the identification map, cf. 2.1;

(b') also,}
$$H^i(C^\ast_A\otimes_P \Z_0) \simeq H^i(N, A\cup \pn; \Z).\tag2-22$$

\subheading{2.5. Algebraic integers and the lifting property} Here it will become clear why our definition of the cup-length $\cl(\xi)$ involves the condition of not being a Dirichlet unit.
\proclaim{Proposition 2} Suppose that $A\subset N$ is a subcomplex such that the inclusion
$A\to N$ is homotopic to a map $A\to \pn$ keeping $A\cap \pn$ fixed. Let $a\in \C^\ast$ be a 
complex number such that $a^{-1}$ is not an algebraic integer. Then the homomorphism
$C^\ast_A \to C^\ast$ induces an epimorphism on the cohomology
$$H^i(C^\ast_A\otimes_P \C_a) \to H^i(C^\ast\otimes_P \C_a), \qquad i=0, 1, 2, \dots.\tag2-23$$ 
\endproclaim
\demo{Proof} Let $\Z_0$ denote the group $\Z$ considered as a $P$-module with the trivial $\tau$ action, i.e. $\Z_0= P/\tau P$. We will show first that
$$H^i(C^\ast_A\otimes_P \Z_0) \to H^i(C^\ast\otimes_P \Z_0)\tag2-24$$
is an epimorphism.
We know from  (b') of subsection 2.4 that
$$H^i(C^\ast_A\otimes_P \Z_0) \simeq H^i(N, A\cup \pn; \Z)\quad\text{and}\quad 
H^i(C^\ast\otimes_P \Z_0) \simeq H^i(N,  \pn;\Z).$$
In the exact sequence
$$\dots \to H^i(N, A\cup \pn;  \Z) \to H^i(N,  \pn;  \Z) 
\overset{j^\ast}\to\longrightarrow 
H^i(A\cup \pn, \pn;\Z)\to \dots $$
$j^\ast$ acts trivially (since the inclusion $(A\cup \pn,\pn)\to (N,\pn)$ is null-homotopic) and hence
$H^i(N, A\cup \pn; \Z) \to H^i(N,  \pn; \Z)$ is an epimorphism. This proves that (2-24) 
is an epimorphism.
Now, Proposition 2 follows from Lemma 5 below. \qed
\enddemo

\proclaim{Lemma 5} Let $C$ and $D$ be chain complexes of free finitely generated 
$P=\Z[\tau]$-modules and let $f:C \to D$ be a chain map. Suppose that for some $q$ 
the induced map $f_\ast: H_q(C\otimes_P \Z_0) \to H_q(D\otimes_P \Z_0)$
is an epimorphism; here $\Z_0$ is $\Z$ considered with the trivial $P$-action: 
$\Z_0= P/\tau P$ . 
Then for any complex number $a\in \C^\ast$, such that $a^{-1}$ is not an algebraic integer, the homomorphism
$$f_\ast: H_q(C\otimes_P\C_a) \to H_q(D\otimes_P\C_a)\tag2-25$$
is an epimorphism; here $\C_a$ denotes $\C$ with $\tau$ acting as the multiplication by $a$.\endproclaim 

\demo{Proof} Denote by $Z_q(C), Z_q(D)$ the cycles of $C$ and $D$ and by $B_q(C)$ and $B_q(D)$ their boundaries. Recall that the homological dimension of $P$ is 2.
We have the exact sequence
$$0\to Z_q(C)\to C_q\to B_{q-1}(C)\to 0$$
and hence $Z_q(C)$ is a free $P$-module (since  $B_{q-1}(C)$ is a submodule of a free module and so has a homological dimension $\le 1$). Similarly $Z_q(D)$ is free. 

Choose a basis for $Z_q(C), Z_q(D)$ and $D_{q+1}$ and express in terms of these basis the map
$$f\oplus d: Z_q(C)\oplus D_{q+1} \to Z_q(D).\tag2-26$$
The resulting matrix $M$ is a rectangular matrix with entries in $P$. 

We claim: {\it there exist  integers $b_j\in \Z$ and  minors $A_j(\tau) \in P$ of the matrix $M$ of size $\rk Z_q(D)\times \rk Z_q(D)$,
such that the polynomial with integer coefficients
$$p(\tau)= \sum_j b_j A_j(\tau)\tag2-27$$ satisfies} 
$$p(0)=1.\tag2-28$$

In fact, we will show that our claim is {\it equivalent} to the requirement that
$f_\ast: H_q(C\otimes_P \Z_0) \to H_q(D\otimes_P \Z_0)$ is an isomorphism.
Namely, using the resolvent 
$0\to P\overset{\tau}\to\longrightarrow P\to \Z_0\to 0$ it is easy to see that 
$\Tor_1^P(B_{q-1}(C),\Z_0)=0$ (since $B_{q-1}(C)$ is a submodule of a free module). Hence we have the exact sequence
$$0\to Z_{q}(C)\otimes_P \Z_0\to C_q\otimes_P\Z_0\to B_{q-1}(C)\otimes \Z_0\to 0.$$
This means that $Z_q(C)\otimes_P \Z_0 = Z_q(C\otimes_P \Z_0)$, and 
$B_{q-1}(C)\otimes_P \Z_0 = B_{q-1}(C\otimes_P \Z_0)$. Hence, the hypothesis of the lemma
implies that 
the homomorphism
$$f\oplus d: (Z_q(C)\otimes_P \Z_0)\oplus (D_{q+1}\otimes_P \Z_0) \to Z_q(D)\otimes_P \Z_0$$
is an epimorphism. This epimorphism is described by the matrix $M(0)$, where we substitute $\tau=0$
into $M$.
Therefore, there are minors $A_j(\tau)$ of $M$ of size $\rk Z_q(D)\times \rk Z_q(D)$ so that the ideal
in $\Z$ generated by the integers $A_j(0)$ contains 1. This proves (2-28). 

Since $p(\tau)$ is an integral polynomial with $p(0)=1$ and $a^{-1}$ is not an algebraic 
integer it follows that 
$$p(a)\ne 0.\tag2-29$$

Let us show that (2-29) is equivalent to the statement that (2-25) is an epimorphism. We have the exact 
sequence
$$0\to Z_{q}(C)\otimes_P \C_a\to C_q\otimes_P\C_a\to B_{q-1}\otimes \C_a\to 0$$
(here we may work over $\C[\tau]$ which is a PID). Hence, similarly to the arguments above, 
we obtain that the map
$$f\oplus d: (Z_q(C)\otimes_P \C_a)\oplus (D_{q+1}\otimes_P \C_a )\to Z_q(D)\otimes_P \C_a\tag2-30$$
is described by the matrix $M$ with substitution $\tau=a$. We conclude that at least one of the 
$\rk Z_q(D)\times \rk Z_q(D)$ minors $A_j(a)$ is nonzero because of (2-29), and hence (2-30) 
and (2-25) are epimorphisms. \qed
\enddemo

\proclaim{2.6. Corollary} Let $a\in \C^\ast$ be a complex number, not an algebraic integer. 
Let $A\subset X$ be a closed subset such that $A=p(A')$, where $A'\subset N-\pn $ is a closed polyhedral subset such that the inclusion $A'\to N$ is homotopic to a map with values in 
$\pn$. Then the 
restriction map
$$H^q(X,A;\E_a)\to H^q(X;\E_a)\tag2-31$$
is an epimorphism.\endproclaim
\demo{Proof} We just combine the isomorphisms (a) and (a') (cf. 2.3, 2.4) and Proposition 2.
\qed \enddemo

\subheading{2.7. End of proof of Theorem 1} We need to establish inequality (2-11). 
In other words,
we want to prove the triviality of any cup-product
$$v_0\cup v_1\cup  \dots \cup v_{m+1} =0,\quad\text{where}\quad v_j\in H^{d_j}(X;\E_{a_j}),
\quad d_j >0,\tag2-32$$
assuming that among the numbers $a_0, a_1, \dots, a_{m+1}\in \C$ at least two are not Dirichlet units; 
here $m$ denotes $m= \cn N)$.

We shall assume that $a_0$ and $a_{m+1}$ are not Dirichlet units; if not, we just rename the numbers.

Moreover, we will assume that one of the numbers 
$a_0$ and $a_{m+1}$ is not an algebraic integer. In the 
case when both $a_0$ and $a_{m+1}$ are algebraic integers the inverse numbers $a_0^{-1}$ 
and $a_{m+1}^{-1}$are not algebraic integers and we shall apply the arguments following below
to the form $-\omega$ (representing the cohomology class $-\xi$), which obviously has the same
set of critical points.) 

Since we may always rename the numbers $a_0$ and $a_{m+1}$, we will assume below
that $a_0$
is not an algebraic integer.

Suppose that $N$ can be covered by closed subsets $A_0, A_1\cup  \dots\cup A_m=N$
so that $A_0$ is a collar of $\pn$ (cf. 2.2), and for $j=1,2, \dots, m$ the subset $A_j$ is disjoint from
$\pn$ and null-homotopic in $N$. Using Lemma 3, we may assume that the sets $A_j$ are 
polyhedral.  We find (since $d_j>0$) that for $j=1,2, \dots, m$ we may lift the class
$v_j$ to a relative cohomology class $\tilde v_j\in H^{d_j}(X,B_j;\E_{a_j})$, where 
$B_j=p(A_j )$.

Let $U_\pm$ be a small cylindrical neighborhood of 
$\partial_\pm N$ in $N$.
Applying Corollary 2.6, class $v_0$ can be lifted to a class $\tilde v_0\in H^{d_j}(X,B_0;\E_{a_0})$,
where $B_0=p(A_0 -U_+)$. 

Let $B_{m+1}$ be a closed cylindrical neighborhood of $V$ in $X$
containing $\overline {p(U_-)}\cup \overline{ p(U_+)}$. We claim that we may lift the class
$v_{m+1}\in H^{d_{m+1}}(X; \E_{a_{m+1}})$ to a class 
$\tilde v_{m+1}\in H^{d_{m+1}}(X, B_{m+1}; \E_{a_{m+1}})$.
We will use Corollary 2.6.
First, find two shifts of $V$ into $X-B_{m+1}$, one (denoted $V'$) in the positive normal direction and the other (denoted $V''$) in the negative normal direction (cf. Figure 5). If the number $a_{m+1}$
is not an algebraic integer we may apply Corollary 2.6 to the cut $V''$. If the number 
$a_{m+1}^{-1}$
is not an algebraic integer we may apply Corollary 2.6 to the cut $V'$.

$$\epsffile{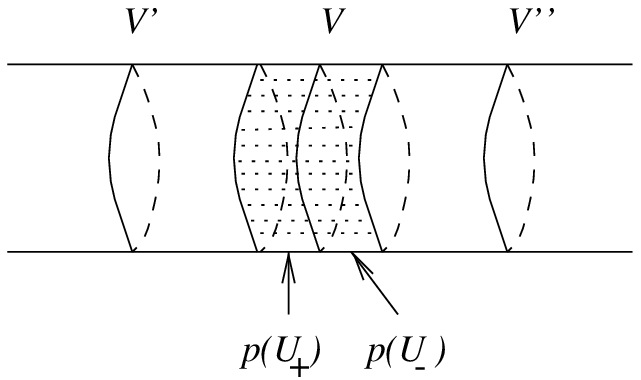}$$
\centerline{Figure 5.}
\smallskip
Now, it is clear that the product 
$v_0\cup\dots \cup v_{m+1}$ is trivial since it is obtained from the product 
$\tilde v_0\cup\dots \cup \tilde v_{m+1}$ (lying in $H^d(X, \cup_{j=0}^{m+1} B_j;\E_a)$, 
where $a=a_0 a_1\dots a_{m+1}$) by restricting onto $X$, and the group $H^d(X, \cup_{j=0}^{m+1}  B_j;\E_a)$
vanishes since $X=\cup_{j=0}^{m+1} B_j$. \qed

\heading{\bf \S 3. Proofs of Theorems 2 and 3 and Proposition 1}\endheading

\subheading{3.1. Proof of Theorem 2} Consider the deformation complex constructed in subsection
2.3. It is a free finitely generated cochain complex $C=\oplus C^q$ over the ring $P=\Z[\tau]$, which
satisfies (2-15). We denote by $B^q(a)$ the rank of the following linear map
$$\delta\otimes 1: C^{q-1}\otimes_P \C_a \to C^{q}\otimes_P \C_a.\tag3-1$$ 
It is clear that this map $\delta\otimes 1$ can be represented by a square matrix whose entries are
polynomials in $a$ with integral coefficients. Therefore we obtain that {\it $B^q(a)$ as a function 
of $a$ is constant except
possibly at finitely many jump points $a= a_1^q, \dots, a_k^q$, which are all algebraic numbers,
and at those jump values $a= a_1^q, \dots, a_k^q$ the value of $B^q(a)$ is smaller than at a generic 
point. } Applying the Euler - Poincar\'e formula to the truncated complex 
$0\to C^q\to C^{q+1}\to\dots$ and also (2-15) we find
$$\aligned
 \dim H^q(X;\E_{a^{-1}}) = &\sum_{i=0}^\infty (-1)^i \rk C^{q+i} -\\
&-\sum_{j=1}^\infty (-1)^j\dim H^{q+j}(X;\E_{a^{-1}}) - B^q(a).
\endaligned\tag3-2
$$ 
The sums here are actually finite. 
This proves statements (1), (3), (4) of Theorem 2. Statement (2) follows from 
the definition of the Novikov numbers given in \cite{F}, subsection 1.2. According to this definition,
the number $b_i(\xi)$ is the dimension over the field of rational functions $F(\tau)$ of the local system of rank 1 over $X$ which is naturally determined by the class $\xi$, cf. \cite{F}, section 1.2.
Therefore, it is obvious that $b_i(\xi)$ coincides with the value $\dim H^i(X,\E_a)$ for a generic point
$a\in \C^\ast$. In section 1.4 of \cite{F}
it is shown that this definition of $b_i(\xi)$ (using the rational functions) 
is in fact equivalent to the original definition 
of S.P. Novikov using the formal power series.\qed

\subheading{3.2. Proof of Theorem 3} Let 
$\omega$ be a closed 1-form lying in a cohomology 
class $\xi\in H^1(X;\R)$ of $\rank =r >1$. Let $S=S(\omega)$ denote the set of zeros of $\omega$. 
It is clear that $\xi|_S =0$.

Let $\xi_1, \dots, \xi_r\in H^1(X;\Z)$ be a basis of the free abelian group 
$\Hom(H_1(X)/\ker (\xi))$, where $r$ is the rank of $\xi$. We may write 
$\xi = \sum_{i=1}^r \alpha_i\xi_i$, and the coefficients are real $\alpha_i\in \R$. 

Suppose that $\xi_m$ is a sequence of rank 1 classes with $\cl(\xi_m)\ge \cl(\xi)$,
which converges to $\xi$ as $m\to \infty$, 
and each of the classes $\xi_m$ vanishes on $\ker(\xi)$. Then we have $\xi_m = \sum_i \alpha_{i,m}\xi_i$, where $\alpha_{i,m}=\lambda_m\cdot n_{i,m}$, $\lambda_m\in \R$, and 
$n_{i,m}\in \Z$ for $i=1, 2, \dots, r$. Each sequence $\alpha_{i,m}$ converges to $\alpha_i$ as
$m$ tends to $\infty$.

Choose a closed 1-form $\omega_i$ in the class $\xi_i$ for $i=1, \dots, r$; since $\xi_i|_S=0$
we may choose it so that it vanishes identically on a neighborhood of $S$. 
Define the following sequence of closed 1-forms
$$\omega_m = \omega - \sum_{i=1}^r (\alpha_i -\alpha_{i,m})\omega_i.$$
It is clear that $\omega_m$ has rank 1 and for $m$ large enough $S(\omega_m)= S(\omega)$. 
The cohomology class of $\omega_m$ is $\xi_m$. By Theorem 1 we have
$\cat(S(\omega)) \ge \cl(\xi_m)-1.$ 
Hence we obtain $\cat(S(\omega)) \ge \cl(\xi)-1$. \qed

\subheading{3.3. Proof of Proposition 1} 
It is clear that it is enough to prove (1-10) assuming that the classes $\xi_1$ and $\xi_2$ are integral
$\xi_\nu\in H^1(X_\nu;\Z)$ for $\nu=1, 2$. The general statement then follows automatically due to the
nature of our definition of $\cl(\xi)$ for general $\xi$, cf. 1.6.

Position $X_1$ and $X_2$ so that their intersection is a small 
$n$-dimensional disk $D^n$, and 
then the connected sum $X_1\# X_2$ is obtained from the union $X_1\cup X_2$ by
removing the interior of $D^n$.
Let $\E$ be a flat bundle over the connected sum
$X_1\# X_2$ and let $\E_\nu$ be a flat bundle over $X_\nu$ so that 
$$\E|_{X_\nu - \overset\circ\to{D^n}}\simeq \E_\nu|_{X_\nu - \overset\circ\to{D^n}},\tag3-3$$ 
for $\nu=1, 2$.
As follows from the Mayer - Vietoris sequence there is a canonical isomorphism
$$\psi:H^q(X_1;\E_1)\oplus H^q(X_2;\E_2)\to H^q(X_1\# X_2;\E)\tag3-4$$
for $0<q<n = \dim X_1=\dim X_2$. For $q=n$ the homomorphism $\psi$ is only 
an epimorphism, but
its restriction on any of the summands $H^n(X_1;\E|_{X_1})$ or $H^n(X_2;\E|_{X_2})$ is an isomorphism. 

Also, $\psi$ is multiplicative in the following sense.
Suppose that we have another flat bundle $\F$ over the connected sum
$X_1\# X_2$ and let $\F_\nu$ be flat bundles over $X_\nu$, $\nu = 1, 2$ satisfying  condition (3-3). Then for any $v\in H^\ast(X_1,\E_1)$ and  $w\in H^\ast(X_1,\F_1)$
holds $\psi(v\cup w, 0)= \psi(v,0)\cup \psi(w,0)$ and similarly with respect to the other variable.

Given a complex number $a\in \C^\ast$, it determines (as in 1.2)
the flat line
bundles $\E_\nu$ over $X_\nu$ (together with the class $\xi_\nu$, $\nu=1, 2$) and a flat line bundle
$\E$ over $X=X_1\# X_2$ (together with $\xi= \xi_1\# \xi_2$). These three
flat bundles clearly satisfy (3-3).

Suppose now that we have classes 
$v_j \in H^{d_j}(X_1;\E_{a_j})$, where $j= 1, 2, \dots, k$ such that
their product $v_1\cup \dots \cup v_k$ is non-trivial. We assume that at least two $a_j$ are not Dirichlet units, and this implies that $d_j<n$ for all $j$. Then we obtain classes $w_j=\psi(v_j, 0)
\in H^{d_j}(X;\E_{a_j})$ with $w_1\cup \dots \cup w_k\ne 0$. 
This proves inequality $\cl(\xi)\ge \cl(\xi_1)$. Therefore $\cl(\xi)\ge \max\{\cl(\xi_1), \cl(\xi_2)\}.$

The inverse inequality follows similarly, using the properties of the map $\psi$ mentioned above.\qed

\subheading{Acknowledgment} It is my pleasure to thank Yuli Rudyak for useful and enjoyable discussions.

\enddocument